\documentclass[11pt,leqno]{amsart}

\usepackage[margin=1.3in]{geometry}

\usepackage[T1]{fontenc}
\usepackage{lmodern}
\usepackage{amsmath,amssymb,amsthm,mathtools,mathrsfs}
\usepackage{enumitem}
\usepackage{hyperref}

\usepackage{graphicx}
\usepackage{amscd}
\usepackage[all]{xy}
\usepackage{tikz}
\usepackage{xcolor}

\newcommand{\F}{\mathbb{F}}
\newcommand{\Q}{\mathbb{Q}}

\newcommand{\Z}{\mathbb{Z}}

\DeclareMathOperator{\Gal}{Gal}
\DeclareMathOperator{\Hom}{Hom}
\DeclareMathOperator{\End}{End}

\DeclareMathOperator{\tr}{tr}

\newcommand{\GL}{\mathrm{GL}}
\newcommand{\Br}{\mathrm{Br}}

\newcommand{\ad}{\mathrm{ad}\,\bar{\rho}}
\newcommand{\adzero}{\mathrm{ad}^{0}\,\bar{\rho}}

\newcommand{\mm}{\mathfrak{m}}

\newcommand{\T}{\mathbb{T}}

\newcommand{\resrho}{\bar{\rho}}
\newcommand{\eps}{\varepsilon}
\newcommand{\epsl}{\varepsilon_\ell}
\newcommand{\barepsl}{\bar{\varepsilon}_\ell}

\DeclareFontFamily{U}{wncy}{}
\DeclareFontShape{U}{wncy}{m}{n}{<->wncyr10}{}
\DeclareSymbolFont{mcy}{U}{wncy}{m}{n}
\DeclareMathSymbol{\Sh}{\mathord}{mcy}{"58}

\theoremstyle{plain}
\newtheorem{theorem}{Theorem}[section]
\newtheorem{proposition}[theorem]{Proposition}
\newtheorem{corollary}[theorem]{Corollary}
\newtheorem{lemma}[theorem]{Lemma}

\theoremstyle{definition}
\newtheorem{definition}[theorem]{Definition}

\newtheorem{assumption}[theorem]{Assumption}

\theoremstyle{remark}
\newtheorem{remark}[theorem]{Remark}

\setlist[itemize]{leftmargin=*,topsep=4pt,itemsep=2pt}
\setlist[enumerate]{leftmargin=*,topsep=4pt,itemsep=2pt}

\title{A Local--Global Study of Obstructed Deformation Problems I}
\author{Bartu Bingol}
\date{}

\begin{document}
\maketitle

\begin{abstract}
We study obstructed deformation problems for two-dimensional residual Galois representations arising from weight~$2$ newforms of level~$N$.
Using Poitou-Tate duality, we isolate local and global sources of obstructions and give concrete criteria for when they occur.
In several cases we also describe the resulting universal deformation ring explicitly.
\end{abstract}


\section{Introduction}

\subsection{Notation and standing assumptions}\label{subsec:notation}

Let $S$ be a finite set of places of $\Q$ containing $\infty$ and all primes dividing $N\ell$.
Let $\Q_S$ be the maximal extension of $\Q$ unramified outside $S$, and set
\[
G_{\Q,S} := \Gal(\Q_S/\Q).
\]
For each prime $p$, we write $G_p \subseteq G_{\Q,S}$ for a decomposition group at $p$ and $I_p\subseteq G_p$ for inertia.

Fix a rational prime $\ell$ and let $\kappa$ be a finite field of characteristic $\ell$, with ring of Witt vectors $W(\kappa)$.
Let $f$ be a normalized Hecke eigenform of weight~$2$ and level~$N$, and let
\[
\rho_f : G_{\Q,S} \longrightarrow \GL_2\bigl(W(\kappa)\bigr)
\]
be an $\ell$-adic Galois representation attached to $f$.
Reducing modulo the maximal ideal of $W(\kappa)$ gives the residual representation
\[
\resrho : G_{\Q,S} \longrightarrow \GL_2(\kappa).
\]
Throughout the paper we assume that $\resrho$ is odd, continuous, and absolutely irreducible.

\subsection{Deformations and deformation rings}\label{subsec:deformations}

Let $\mathcal{C}$ denote the category of complete Noetherian local $W(\kappa)$-algebras $(R,\mm_R)$ with residue field $R/\mm_R \cong \kappa$.

\begin{definition}\label{def:lifts-deformations}
Let $R\in\mathcal{C}$.
\begin{enumerate}
    \item A \emph{lift} of $\resrho$ to $R$ is a continuous homomorphism
    \[
    \rho_R: G_{\Q,S}\to \GL_2(R)
    \]
    whose reduction modulo $\mm_R$ equals $\resrho$.
    \item Two lifts $\rho_R$ and $\rho_R'$ are \emph{strictly equivalent} if there exists $g\in 1+M_2(\mm_R)$ such that
    $\rho_R' = g\rho_R g^{-1}$.
    \item A \emph{deformation} of $\resrho$ to $R$ is a strict equivalence class of lifts.
\end{enumerate}
\end{definition}

We write $D_{\resrho}(R)$ for the set of deformations of $\resrho$ to $R$.
Under the standard hypotheses of Mazur~\cite{Maz87}, the functor
\[
D_{\resrho}:\mathcal{C}\to \mathbf{Sets}, \qquad R\mapsto D_{\resrho}(R)
\]
is pro-representable by a complete local ring $R^{\mathrm{univ}}$ together with a universal deformation
\[
\rho^{\mathrm{univ}}:G_{\Q,S}\to \GL_2(R^{\mathrm{univ}}).
\]
Equivalently, for each $R\in\mathcal{C}$ there is a natural identification
\[
D_{\resrho}(R)\ \cong\ \Hom_{\mathcal{C}}\!\left(R^{\mathrm{univ}},R\right).
\]

Congruences between eigenforms are naturally encoded by Hecke algebras. Let $\T$ denote the Hecke algebra acting on weight~$2$ cusp forms of level~$N$, and let $\mm\subseteq \T$ be the maximal ideal cut out by the residual eigenvalues of $f$ modulo~$\ell$.
Eigenforms congruent to $f$ modulo~$\ell$ correspond to the same maximal ideal $\mm$, and $\mm$ controls $\resrho$.

\subsection{Outline and main results}\label{subsec:outline}

Our main structural tool is the short exact sequence extracted from Poitou-Tate duality:
\begin{equation}\label{eq:PT}
0 \longrightarrow \prod_{p\in S} H^0\!\left(G_p,\ \barepsl\otimes \ad(\resrho)\right)
\longrightarrow \Hom_\kappa\!\left(H^2(G_{\Q,S},\ad(\resrho)),\kappa\right)
\longrightarrow \Sh^1\!\left(G_{\Q,S},\barepsl\otimes \ad(\resrho)\right)
\longrightarrow 0.
\end{equation}
Here $\barepsl$ denotes the mod~$\ell$ cyclotomic character and $\Sh^1$ is the usual Tate-Shafarevich-style kernel.

In \S\ref{sec:deformation-theory} we recall the deformation functor, its tangent space, and the obstruction group $H^2(G_{\Q,S},\ad(\resrho))$.
In \S\ref{sec:level-raising} we study a clean level-raising mechanism: allowing ramification at a new prime $q\notin S$ creates a \emph{single} local obstruction, and we compute it explicitly (including a concrete cocycle formula via a Brauer-group recipe).
Finally, in \S\ref{sec:global-obstructions} we turn to genuinely global obstructions coming from $\Sh^1$ and explain how they are governed by strict congruences between newforms.

A recurring theme is the difference between \emph{optimal} and \emph{non-optimal} level. Roughly speaking, if a residual representation arises from level~$N$ but actually has a congruent newform at a smaller level, then one should expect obstructions to show up. The precise level-lowering criteria that detect when a mod-$\ell$ representation has non-optimal level are developed in Diamond--Taylor~\cite{DT94}, and Hatley~\cite{Hat16} shows that non-optimal level forces obstructed modular deformation problems under standard hypotheses. We keep these results in the background as guiding principles, and we point to them explicitly in \S\ref{sec:global-obstructions} where ``optimal level'' becomes part of the obstruction story.

We now state one of our main level-raising results; the explicit cocycle formula is given in \S\ref{subsec:diagram-recipe}.

\begin{theorem}\label{thm:levelRaising}
Let $\resrho:G_{\Q,S}\to \GL_2(\F_\ell)$ arise from a weight~$2$ newform of level~$N$, where $S$ contains the primes dividing $N\ell$, and assume $\ell>3$.
Suppose the deformation problem for $\resrho$ over $G_{\Q,S}$ is unobstructed.
Let $q\notin S$ be a prime such that $q^2\not\equiv 1\pmod{\ell}$, and enlarge $S' := S\cup\{q\}$.

Then the deformation problem for $\resrho$ over $G_{\Q,S'}$ is obstructed. More precisely:
\begin{enumerate}
\item The obstruction group $H^2(G_{\Q,S'},\ad(\resrho))$ is one-dimensional over $\F_\ell$.
\item A generator may be chosen so that the obstruction is detected by a single matrix entry (the upper-right entry in a convenient basis). Concretely, this is governed by an explicit exponent formula in tame generators (see \S\ref{subsec:diagram-recipe}).
\item The universal deformation ring for $G_{\Q,S'}$ is of the form
\[
\Z_\ell[[T_1,T_2,T_3,T]]/\langle T(\ell-\Phi)\rangle,
\]
where $\Phi$ is an explicit Frobenius power series determined by the local computation at $q$.
\end{enumerate}
\end{theorem}

We also obtain results describing when global obstructions occur in the ``optimal level'' setting.
In particular, when $f$ has rational Fourier coefficients (equivalently, $f$ is attached to an elliptic curve $E_f/\Q$), strict congruence primes for $f$ play a key role in controlling $\Sh^1$.


\section{Deformation theory of two-dimensional Galois representations}\label{sec:deformation-theory}

\subsection{The deformation functor and its tangent space}\label{subsec:tangent}

For the remainder of the paper, unless stated otherwise, we fix $\resrho:G_{\Q,S}\to \GL_2(\kappa)$ as in \S\ref{subsec:notation}.
Our first goal is to connect infinitesimal deformations of $\resrho$ to Galois cohomology.

Let $\kappa[\eps]/(\eps^2)$ be the ring of dual numbers, and consider a lift
\[
\rho_\eps: G_{\Q,S}\to \GL_2\bigl(\kappa[\eps]/(\eps^2)\bigr)
\]
reducing to $\resrho$ modulo $\eps$.
Writing $\rho_\eps(\sigma)=(I+\eps\,b(\sigma))\resrho(\sigma)$ for a function $b:G_{\Q,S}\to M_2(\kappa)$, a standard computation shows that $\rho_\eps$ is a homomorphism if and only if $b$ is a $1$-cocycle for the conjugation action of $G_{\Q,S}$ on $M_2(\kappa)$ via $\resrho$.
This $G_{\Q,S}$-module is the adjoint representation $\ad(\resrho)$.

\begin{proposition}\label{prop:tangentH1}
There is a natural isomorphism of $\kappa$-vector spaces
\[
D_{\resrho}(\kappa[\eps]/(\eps^2))\ \cong\ H^1(G_{\Q,S},\ad(\resrho)).
\]
\end{proposition}

\begin{proof}
A lift $\rho_\eps$ determines a function $b$ as above, and the homomorphism condition for $\rho_\eps$ is equivalent to the cocycle condition for $b$.
Changing $\rho_\eps$ by strict equivalence changes $b$ by a coboundary.
Thus strict equivalence classes of lifts to $\kappa[\eps]/(\eps^2)$ are identified with $H^1(G_{\Q,S},\ad(\resrho))$.
\end{proof}

As a useful consequence we recover the standard ``power series upper bound'' for $R^{\mathrm{univ}}$ (cf.\ \cite[\S2.6]{DDT95}, \cite{Gou01}).

\begin{corollary}\label{cor:power-series}
Assume $R^{\mathrm{univ}}$ exists.
Let $d_1:=\dim_\kappa H^1(G_{\Q,S},\ad(\resrho))$.
Then $R^{\mathrm{univ}}$ is a quotient of $W(\kappa)[[T_1,\dots,T_{d_1}]]$.
\end{corollary}

\subsection{Obstructions}\label{subsec:obstructions}

Let $\psi:R_2\twoheadrightarrow R_1$ be a surjection in $\mathcal{C}$ with kernel $\mathfrak{a}$ satisfying $\mathfrak{a}\mm_{R_2}=0$.
Fix a deformation $\rho_{R_1}:G_{\Q,S}\to \GL_2(R_1)$, and choose a set-theoretic lift $\widetilde{\rho}:G_{\Q,S}\to \GL_2(R_2)$ reducing to $\rho_{R_1}$.

The defect of $\widetilde{\rho}$ being a homomorphism is measured by
\[
c_{\widetilde{\rho}}:G_{\Q,S}\times G_{\Q,S}\longrightarrow 1+M_2(\mathfrak{a}),
\qquad
c_{\widetilde{\rho}}(\sigma,\tau)
=\widetilde{\rho}(\sigma\tau)\widetilde{\rho}(\sigma)^{-1}\widetilde{\rho}(\tau)^{-1}.
\]
Writing $c_{\widetilde{\rho}}(\sigma,\tau)=1+d_{\widetilde{\rho}}(\sigma,\tau)$ identifies
$d_{\widetilde{\rho}}$ with a $2$-cochain valued in $\ad(\resrho)\otimes_\kappa \mathfrak{a}$.
A standard computation shows that $d_{\widetilde{\rho}}$ is a $2$-cocycle, and its cohomology class is the obstruction to lifting $\rho_{R_1}$ across $\psi$.
In particular, obstructions lie in $H^2(G_{\Q,S},\ad(\resrho))$.

\begin{definition}
The deformation problem for $\resrho$ over $G_{\Q,S}$ is \emph{unobstructed} if $H^2(G_{\Q,S},\ad(\resrho))=0$, and \emph{obstructed} otherwise.
\end{definition}

A useful unobstructedness criterion in the weight~$2$ (elliptic curve) setting is due to Flach~\cite{Fla92}.
For complementary obstruction criteria tailored to modular deformation problems (including several non-optimal-level scenarios),
see Hatley~\cite{Hat16}.


\section{Obstructions arising from level raising}\label{sec:level-raising}

\subsection{Global setup and the Poitou--Tate reduction}\label{subsec:S-to-Sq}

Fix a finite set of places $S$ of $\Q$ containing $\infty$ and all primes dividing $N\ell$.
Let $\Q_S$ be the maximal extension of $\Q$ unramified outside $S$, and set
\[
G_{\Q,S}:=\Gal(\Q_S/\Q).
\]
For each prime $p$ we write $G_p\subseteq G_{\Q,S}$ for a decomposition group at $p$.
Let $\barepsl:G_{\Q,S}\to \F_\ell^\times$ be the mod-$\ell$ cyclotomic character.

Let
\[
\bar{\rho}:G_{\Q,S}\longrightarrow \GL_2(\F_\ell)
\]
be odd, continuous, and absolutely irreducible, arising from a weight $2$ newform of level $N$.
Write
\[
\ad := \End_{\F_\ell}(\F_\ell^2)
\]
with $G_{\Q,S}$-action by conjugation via $\bar{\rho}$, and let $\adzero\subset \ad$
be the trace-zero subrepresentation.

We now enlarge the ramification set by a prime $q\notin S$. Put
\[
S' := S\cup\{q\}, \qquad G_{\Q,S'}:=\Gal(\Q_{S'}/\Q),
\]
where $\Q_{S'}$ is the maximal extension unramified outside $S'$.
The obstruction space for the deformation problem over $G_{\Q,S'}$ is
\[
H^2\!\left(G_{\Q,S'},\ad\right).
\]

The key input is the following short exact sequence extracted from Poitou--Tate duality
(see \cite[Thm.~4.10]{Mil86} for a full statement):
\begin{equation}\label{eq:PT-sec3}
0 \longrightarrow \prod_{p\in S'}
H^0\!\left(G_p,\barepsl\otimes \ad\right)
\longrightarrow
\Hom_{\F_\ell}\!\left(H^2(G_{\Q,S'},\ad),\F_\ell\right)
\longrightarrow
\Sh^1\!\left(G_{\Q,S'},\barepsl\otimes \ad\right)
\longrightarrow 0,
\end{equation}
where
\[
\Sh^1(G_{\Q,S'},M):=\ker\!\left(H^1(G_{\Q,S'},M)\to \prod_{p\in S'}H^1(G_p,M)\right).
\]

In this section we treat a clean level-raising situation in which the obstruction over $S'$ comes \emph{only} from the new local term at $q$.
Accordingly we impose:

\begin{assumption}\label{assump:unobstructed-over-S}
The deformation problem for $\bar{\rho}$ over $G_{\Q,S}$ is unobstructed, i.e.
\[
H^2(G_{\Q,S},\ad)=0.
\]
Equivalently (by \eqref{eq:PT-sec3} with $S$ in place of $S'$),
\[
\prod_{p\in S}H^0\!\left(G_p,\barepsl\otimes \ad\right)=0
\quad\text{and}\quad
\Sh^1\!\left(G_{\Q,S},\barepsl\otimes \ad\right)=0.
\]
\end{assumption}

With Assumption~\ref{assump:unobstructed-over-S}, the only potential new contribution to the left term of \eqref{eq:PT-sec3}
is at $p=q$. Thus \eqref{eq:PT-sec3} specializes to
\begin{equation}\label{eq:PT-reduced-sec3}
0 \longrightarrow
H^0\!\left(G_q,\barepsl\otimes \ad\right)
\longrightarrow
\Hom_{\F_\ell}\!\left(H^2(G_{\Q,S'},\ad),\F_\ell\right)
\longrightarrow
\Sh^1\!\left(G_{\Q,S'},\barepsl\otimes \ad\right)
\longrightarrow 0.
\end{equation}

In the level-raising situation we focus on, $\Sh^1(G_{\Q,S'},\barepsl\otimes \ad)=0$
and the obstruction is detected by the injection from the local invariant group at $q$.
So our main tasks are:
\begin{enumerate}[label=\textup{(\arabic*)}, leftmargin=2.2em]
\item compute $H^0(G_q,\barepsl\otimes \adzero)$ explicitly and show it is one-dimensional;
\item explain how this local invariant canonically pairs with the obstruction space and produces a concrete generator;
\item compute that generator in tame generators and obtain an explicit exponent formula.
\end{enumerate}

\subsection{Tame generators at $q$ and the invariant line}\label{subsec:invariants-at-q}

Assume $\ell>3$ and let $q\notin S$ be a prime such that
\begin{equation}\label{eq:q-condition}
q^2\not\equiv 1\pmod{\ell}.
\end{equation}
We also assume $\bar{\rho}$ is unramified at $q$ and that $q$ satisfies the usual level-raising eigenvalue condition
for weight $2$ (so that the ratio of Frobenius eigenvalues matches $\barepsl(F_q)$).
Concretely, after choosing a basis we may suppose
\begin{equation}\label{eq:Frob-diag}
\bar{\rho}(F)=
\begin{pmatrix}\alpha&0\\0&\beta\end{pmatrix},
\qquad
\frac{\alpha}{\beta}\equiv \barepsl(F)\equiv q \pmod{\ell},
\end{equation}
where $F\in G_q$ is an arithmetic Frobenius element.

Since $\bar{\rho}$ is unramified at $q$, inertia acts trivially on $\adzero$.
We now fix tame generators as follows.
Let $I_q\subseteq G_q$ be inertia and let $\tau$ denote a tame inertia element of order $\ell$ in a tame quotient.
Then in the tame quotient we have the standard relation
\begin{equation}\label{eq:tame-relation}
F\tau F^{-1}=\tau^{q}.
\end{equation}
Every element of the tame quotient can be written uniquely as $F^i\tau^{i'}$ with $i\in \Z$ and $i'\in \Z/\ell\Z$.

Choose the standard basis of trace-zero matrices
\[
e_1=
\begin{pmatrix}1&0\\0&-1\end{pmatrix},
\quad
e_2=
\begin{pmatrix}0&1\\0&0\end{pmatrix},
\quad
e_3=
\begin{pmatrix}0&0\\1&0\end{pmatrix}.
\]

\begin{lemma}\label{lem:H0-dim1}
Under \eqref{eq:q-condition} and \eqref{eq:Frob-diag},
\[
\dim_{\F_\ell} H^0\!\left(G_q,\barepsl\otimes \adzero\right)=1,
\]
and the invariant line is generated by $e_3$ (viewed in $\barepsl\otimes \adzero$).
\end{lemma}

\begin{proof}
Because inertia acts trivially on $\adzero$, invariance is determined by the action of $F$.
The $G_q$-action on $\barepsl\otimes \adzero$ is
\[
F\cdot X \;=\; \barepsl(F)\,\bar{\rho}(F)\,X\,\bar{\rho}(F)^{-1}.
\]
From \eqref{eq:Frob-diag}, a direct computation gives
\[
\bar{\rho}(F)e_1\bar{\rho}(F)^{-1}=e_1,\quad
\bar{\rho}(F)e_2\bar{\rho}(F)^{-1}=(\alpha/\beta)e_2,\quad
\bar{\rho}(F)e_3\bar{\rho}(F)^{-1}=(\beta/\alpha)e_3.
\]
Multiplying by $\barepsl(F)\equiv q$ yields
\[
F\cdot e_1 = q\,e_1,\qquad
F\cdot e_2 \equiv q(\alpha/\beta)e_2 \equiv q^2 e_2,\qquad
F\cdot e_3 \equiv q(\beta/\alpha)e_3 \equiv e_3.
\]
Since $q\not\equiv 1\pmod{\ell}$ and $q^2\not\equiv 1\pmod{\ell}$ by \eqref{eq:q-condition}, only $e_3$ is fixed.
\end{proof}

\subsection{Pairings and the obstruction detected by an exponent}\label{subsec:pairing-setup}

The trace pairing $(X,Y)\mapsto \tr(XY)$ identifies $(\adzero)^\vee \simeq \adzero$
as $G_q$-modules. After twisting, local Tate duality provides a perfect pairing
\begin{equation}\label{eq:local-duality-sec3}
H^0\!\left(G_q,\barepsl\otimes \adzero\right)
\times
H^2\!\left(G_q,\adzero\right)
\longrightarrow
H^2(G_q,\mu_\ell)\ \cong\ \Z/\ell\Z.
\end{equation}
Thus $H^2(G_q,\adzero)$ is one-dimensional over $\F_\ell$, and any class in it is determined by its pairing with the generator $e_3$.

To keep the later computation concrete, we package \eqref{eq:local-duality-sec3} as a two-step map:
\begin{equation}\label{eq:cup-trace-map}
H^0\!\left(G_q,\barepsl\otimes \adzero\right)
\ \otimes\
H^2\!\left(G_q,\adzero\right)
\xrightarrow{\ \smile\ }
H^2\!\left(G_q,\barepsl\otimes \adzero\otimes \adzero\right)
\xrightarrow{\ \tr\ }
H^2(G_q,\mu_\ell),
\end{equation}
where $\smile$ is cup product and the last arrow is induced by $(X,Y)\mapsto \tr(XY)$ together with $\barepsl\otimes \F_\ell\simeq \mu_\ell$.

In practice, this means the obstruction class can be read off from a $\mu_\ell$-valued $2$-cocycle on $G_q$,
and we can record it by an exponent function
\[
b(\sigma_1,\sigma_2)\in \F_\ell
\quad\text{such that}\quad
B(\sigma_1,\sigma_2)=\zeta_\ell^{\,b(\sigma_1,\sigma_2)}
\]
for a fixed primitive $\ell$th root of unity $\zeta_\ell$.

\subsection{The diagram and the explicit recipe in the Brauer group}\label{subsec:diagram-recipe}

We now compute the exponent function $b$ explicitly in the tame generators $F,\tau$.
The computation is most transparent when organized through the Brauer group of a local field.

Let $k:=\Q_q$ and $G_k:=\Gal(\overline{k}/k)$.
Recall that
\[
\Br(k):=H^2(G_k,\overline{k}^{\times})\cong \Q/\Z,
\qquad
\Br(k)[\ell]\cong H^2(G_k,\mu_\ell)\cong \Z/\ell\Z.
\]
Let $k^{\mathrm{nr}}$ be the maximal unramified extension of $k$, write
\[
\widetilde{G}_k:=\Gal(k^{\mathrm{nr}}/k)\cong \widehat{\Z},
\]
and let $F\in \widetilde{G}_k$ denote arithmetic Frobenius.

We will use the following standard commutative diagram (see Chapter~XIII of \cite{Ser79} and \cite[Thm.~7.2.6]{NSW08}),
in which the \emph{top row} encodes the invariant map for the Brauer group and the \emph{bottom row} encodes the cup product pairing:

\[
\xymatrixcolsep{2.3pc}\xymatrix{
\Br(k) &
H^2(\widetilde{G}_k,(k^{\mathrm{nr}})^\times) \ar[l]_{\alpha} \ar[r]^{\beta} &
H^2(\widetilde{G}_k,\Z) &
H^1(\widetilde{G}_k,\Q/\Z) \ar[l]_{\delta} \ar[r]^{\gamma} &
\Q/\Z \\
&
H^2(G_k,\mu_\ell) \ar[u]^{\simeq} &
H^2(G_k,M\otimes M^\vee) \ar[l]_{\tr} &
H^1(G_k,M)\otimes H^1(G_k,M^\vee) \ar[l]_-{\smile}
}
\]
Here $M$ is a finite $G_k$-module of $\ell$-power order (in our application, $M=\adzero$), and $M^\vee$ is its Cartier dual. 

\vspace{0.4em}
\noindent\textbf{Conventions and identifications.}
We fix once and for all a primitive $\ell$th root of unity $\zeta_\ell\in \overline{k}$, and we identify
$\mu_\ell=\langle \zeta_\ell\rangle$ with $\F_\ell(1)$ via $\zeta_\ell^{\,a}\leftrightarrow a$.
With this normalization, the mod-$\ell$ cyclotomic character $\bar{\varepsilon}_\ell:G_k\to \F_\ell^\times$
records the action on $\mu_\ell$ by
\[
\sigma(\zeta_\ell)=\zeta_\ell^{\,\bar{\varepsilon}_\ell(\sigma)} \qquad (\sigma\in G_k).
\]
Equivalently, the natural $G_k$-module isomorphism $\F_\ell(1)\simeq \mu_\ell$ gives a canonical identification
\[
\bar{\varepsilon}_\ell\otimes_{\F_\ell}\F_\ell \;\cong\; \mu_\ell,
\]
and throughout we implicitly use this to view $\bar{\varepsilon}_\ell$-twisted scalar pairings as $\mu_\ell$-valued.

\vspace{0.4em}
\noindent\textbf{On generators and the invariant map.}
Under the standard invariant isomorphism $\text{inv}_k:\Br(k)\xrightarrow{\sim}\Q/\Z$,
the class of $1/\ell\in \Q/\Z$ corresponds to a generator of the $\ell$-torsion subgroup
$\Br(k)[\ell]\cong H^2(G_k,\mu_\ell)$.
In the recipe below, the cocycle $C_1$ is chosen so that its image in $\Br(k)$ has invariant $1/\ell$,
hence its class maps to a generator of $H^2(G_k,\mu_\ell)$.

\vspace{0.4em}
\noindent\textbf{The ``divide by a coboundary'' step.}
After passing to $K/k$ with $q=(q^{1/\ell})^\ell\in K^{\times\ell}$, the inflated cocycle $C$ becomes cohomologous
to a $\mu_\ell$-valued cocycle: explicitly, multiplying $C$ by the inverse of a coboundary $\partial\gamma$
does not change its class in $H^2(G,\overline{k}^{\times})$, and because $q$ is an $\ell$th power in $K^\times$,
the resulting cocycle $B:=C/(\partial\gamma)$ takes values in $\mu_\ell$.
Thus $[B]\in H^2(G,\mu_\ell)$ is a well-defined cohomology class, independent of the auxiliary choices up to cohomology,
and it represents the same generator as $[C_1]$ under inflation.

\vspace{0.5em}
\noindent
\textbf{Recipe.}
We now follow the diagram to build an explicit generator of $H^2(G_k,\mu_\ell)$ and then read off the exponent function.

\begin{enumerate}[label=\textup{Step \arabic*:}, leftmargin=2.6em]

\item \textbf{Start from $1/\ell\in \Q/\Z$ and pull back to $H^1(\widetilde{G}_k,\Q/\Z)$.}
Let $f_1\in Z^1(\widetilde{G}_k,\Q/\Z)$ be defined by
\[
f_1(F^i)=\frac{i}{\ell}\ \ (\mathrm{mod}\ \Z).
\]
Thus $\gamma(f_1)=1/\ell$.

\item \textbf{Apply the connecting map $\delta$ to obtain an integral $2$-cocycle.}
Let $\widetilde{f}_1:=\delta(f_1)\in Z^2(\widetilde{G}_k,\Z)$.
Writing $\{x\}=x-\lfloor x\rfloor$ for fractional part, one checks
\[
\widetilde{f}_1(F^i,F^j)=
\left\{\frac{i+j}{\ell}\right\}-\left\{\frac{i}{\ell}\right\}-\left\{\frac{j}{\ell}\right\}.
\]

\item \textbf{Push forward to a $(k^{\mathrm{nr}})^\times$-valued cocycle.}
Let $C_1\in Z^2(\widetilde{G}_k,(k^{\mathrm{nr}})^\times)$ be given by
\[
C_1(F^i,F^j):=q^{\,\widetilde{f}_1(F^i,F^j)}.
\]
Then $[C_1]$ corresponds to $1/\ell\in \Br(k)\cong \Q/\Z$, hence generates $\Br(k)[\ell]$.

\item \textbf{Refine the field so that $q$ becomes an $\ell$th power and isolate the $\mu_\ell$-part.}
Set
\[
K:=k^{\mathrm{nr}}(\zeta_\ell,q^{1/\ell}).
\]
Let $\tau$ be the element with $\tau(q^{1/\ell})=\zeta_\ell q^{1/\ell}$ and $\tau(\zeta_\ell)=\zeta_\ell$,
and extend Frobenius $F$ by $F(q^{1/\ell})=q^{1/\ell}$ and $F(\zeta_\ell)=\zeta_\ell^{q}$.
Then (in the tame quotient) $F\tau F^{-1}=\tau^{q}$.

Inflate $C_1$ to $G:=\Gal(K/k)$ by declaring it trivial on inertia:
\[
C(F^i\tau^{i'},F^j\tau^{j'}) := C_1(F^i,F^j).
\]
Since $q=(q^{1/\ell})^\ell$ is an $\ell$th power in $K^\times$, we can divide $C$ by a coboundary so that the result becomes $\mu_\ell$-valued.

Define a $1$-cochain $\gamma\in C^1(G,K^\times)$ by
\[
\gamma(F^i\tau^{i'}) := (q^{1/\ell})^{\langle i\rangle},
\]
where $\langle i\rangle\in\{0,1,\dots,\ell-1\}$ is the reduction of $i$ modulo $\ell$.
A direct computation using $F\tau F^{-1}=\tau^q$ shows
\[
\frac{C}{\partial\gamma}\!\left(F^i\tau^{i'},F^j\tau^{j'}\right)=\zeta_\ell^{\,i'j q^{i}}.
\]
Thus the $\mu_\ell$-valued cocycle
\[
B(\sigma_1,\sigma_2):=\frac{C(\sigma_1,\sigma_2)}{(\partial\gamma)(\sigma_1,\sigma_2)}
\]
represents a generator of $H^2(G_k,\mu_\ell)\cong \Z/\ell\Z$.

\item \textbf{Read off the exponent function.}
Writing $B(\sigma_1,\sigma_2)=\zeta_\ell^{\,b(\sigma_1,\sigma_2)}$ yields the explicit formula
\begin{equation}\label{eq:explicit-b}
b\!\left(F^i\tau^{i'},\,F^j\tau^{j'}\right)\equiv i'j q^{i}\pmod{\ell},
\qquad (i,j\in\Z,\ i',j'\in\Z/\ell\Z).
\end{equation}

\end{enumerate}

\subsection{Consequences for the obstruction space and the deformation ring}\label{subsec:consequences-sec3}

We now connect the explicit cocycle \eqref{eq:explicit-b} to the obstruction space.

Let $m\in H^0(G_q,\barepsl\otimes \adzero)$ be the generator represented by $e_3$
from Lemma~\ref{lem:H0-dim1}. By the perfect pairing \eqref{eq:local-duality-sec3},
there is a unique (up to scalar) class
\[
[f_q]\in H^2(G_q,\adzero)
\]
whose image under the composite \eqref{eq:cup-trace-map} is the generator of $H^2(G_q,\mu_\ell)$ represented by $B$ above.
In particular, $H^2(G_q,\adzero)$ is cyclic of order $\ell$, and the class $[f_q]$ is detected by the exponent function $b$ in \eqref{eq:explicit-b}.

Finally, returning to the global exact sequence \eqref{eq:PT-reduced-sec3}:
under the level-raising hypotheses (so that the $\Sh^1$-term vanishes and the local invariant at $q$ is nontrivial),
the injection
\[
H^0\!\left(G_q,\barepsl\otimes \ad\right)\hookrightarrow
\Hom_{\F_\ell}\!\left(H^2(G_{\Q,S'},\ad),\F_\ell\right)
\]
forces $H^2(G_{\Q,S'},\ad)$ to be one-dimensional over $\F_\ell$.
Thus the deformation problem over $G_{\Q,S'}$ has a \emph{single} obstruction relation.
The resulting one-relation shape of the universal deformation ring in this $h^2=1$ situation
is treated in detail in \cite{Bos92} and \cite{Boe99}; we record the explicit presentation as part of Theorem~\ref{thm:levelRaising} in the introduction.

\begin{remark}[Local level-raising vs.\ non-optimal level]\label{rem:sec3-nonoptimal}
The obstruction constructed here is genuinely local at a new prime $q$.
This is morally different from the ``non-optimal level'' phenomenon: there, the residual representation arises from level~$N$
but admits a modular realization at smaller level, and one expects obstructions without adding any new prime.
Criteria detecting non-optimal level (especially in squarefree situations) are developed in Diamond--Taylor~\cite{DT94},
and Hatley~\cite{Hat16} proves that non-optimal level forces obstructed modular deformation problems under standard hypotheses. In the next section, we treat several cases regarding this global obstruction mechanism.
\end{remark}


\section{Global obstructions and strict congruence primes}\label{sec:global-obstructions}
\subsection*{Main results of this section}

We isolate the global contribution to obstructions by separating it cleanly from the local invariant terms in the Poitou--Tate exact sequence.

\begin{proposition}[Poitou--Tate reduction to \texorpdfstring{$\Sh^1$}{Sha1}]\label{prop:sec4-main-PT}
Assume \eqref{eq:sec4-standing}. If
\[
H^0\!\left(G_p,\ \bar{\varepsilon}_\ell\otimes \mathrm{ad}\right)=0
\qquad\text{for every }p\in S,
\]
then Poitou--Tate duality yields a natural identification
\[
\Hom_{\F_\ell}\!\left(H^2(G_{\Q,S},\mathrm{ad}(\bar{\rho})),\F_\ell\right)
\ \cong\
\Sh^1\!\left(G_{\Q,S},\bar{\varepsilon}_\ell\otimes \mathrm{ad}(\bar{\rho})\right).
\]
In particular, under the same hypotheses, the deformation problem for $\bar{\rho}$ is obstructed if and only if
\[
\Sh^1\!\left(G_{\Q,S},\bar{\varepsilon}_\ell\otimes \mathrm{ad}(\bar{\rho})\right)\neq 0.
\]
\end{proposition}

\begin{corollary}[Strict congruences force global obstructions]\label{cor:sec4-main-cong}
Assume \eqref{eq:sec4-standing} and the local vanishing condition in Proposition~\ref{prop:sec4-main-PT}.
Under the Selmer hypotheses of Weston~\cite{Wes04,Wes05}, one has
\[
\ell \text{ is a strict congruence prime for } f
\quad\Longrightarrow\quad
\Sh^1\!\left(G_{\Q,S},\bar{\varepsilon}_\ell\otimes \mathrm{ad}(\bar{\rho})\right)\neq 0,
\]
and hence the deformation problem for $\bar{\rho}$ is obstructed.
\end{corollary}

Section~\ref{sec:level-raising} focused on a deliberately clean \emph{local} mechanism: we enlarged the ramification set and watched a new local invariant at a prime $q\notin S$ inject into the obstruction space.
In the present section we turn to the complementary phenomenon.
Here we keep the local terms under control (in particular, we arrange vanishing of the groups
$H^0(G_p,\bar{\varepsilon}_\ell\otimes \mathrm{ad}(\bar{\rho}))$ for $p\in S$),
so that Proposition~\ref{prop:sec4-main-PT} forces any obstruction to be genuinely \emph{global}, detected by
\[
\Sh^1\!\left(G_{\Q,S},\bar{\varepsilon}_\ell\otimes \mathrm{ad}(\bar{\rho})\right).
\]
Once the problem has been reduced to $\Sh^1$, strict congruence primes enter naturally: under the Selmer hypotheses of
\cite{Wes04,Wes05}, Corollary~\ref{cor:sec4-main-cong} explains that strict congruences at level $N$ produce global obstructions.
The rest of the section is then organized around two practical questions:
how to rule out local contributions in concrete situations, and how to translate a global obstruction into explicit deformation-theoretic information.

\subsection{Notation and standing hypotheses}\label{subsec:sec4-setup}

Let $N\ge 1$ and let $\ell$ be a prime.
Let
\[
S:=\{\infty\}\ \cup\ \{p:\ p\mid N\ell\},
\qquad
G_{\Q,S}:=\Gal(\Q_S/\Q),
\]
where $\Q_S$ denotes the maximal extension of $\Q$ unramified outside $S$.
For each prime $p\in S$ let $G_p\subseteq G_{\Q,S}$ be a decomposition group at $p$.

Fix a residual representation
\[
\bar{\rho}:G_{\Q,S}\longrightarrow \GL_2(\F_\ell),
\]
assumed odd, continuous, and absolutely irreducible, arising from a weight $2$ newform $f$ of level $N$.
We again write $\barepsl$ for the mod-$\ell$ cyclotomic character and use
\[
\ad := \End_{\F_\ell}(\F_\ell^2),
\qquad
\adzero := \{X\in \ad : \tr(X)=0\},
\]
with $G_{\Q,S}$-action by conjugation through $\bar{\rho}$.

For the remainder of this section we work under the following mild restrictions, which are standard in the deformation-theoretic references we use:
\begin{equation}\label{eq:sec4-standing}
\ell\nmid N,\qquad \ell\neq 2,\qquad p\neq 2\ \text{and}\ p\not\equiv 1\pmod{\ell}\ \text{for all }p\mid N.
\end{equation}

\subsection{Strict congruence primes and the global obstruction group}\label{subsec:Sha1-congruence}

We begin by defining the congruence condition that governs the global obstruction term.

\begin{definition}[Congruence primes and strict congruence primes]\label{def:sec4-cong}
Let $f$ be a weight $2$ newform of level $N$. A prime $\ell$ is a \emph{congruence prime for $f$}
if there exists a weight $2$ newform $g$ of level $d\mid N$ such that $g$ is not Galois-conjugate to $f$ and
\[
a_n(f)\equiv a_n(g)\pmod{\lambda}
\]
for all but finitely many $n$, for some prime $\lambda$ above $\ell$ in a coefficient field for $f$ and $g$.
If one can take $d=N$, then $\ell$ is called a \emph{strict congruence prime} for $f$.
\end{definition}

The next statement explains why strict congruence primes are exactly the situations in which the global term
$\Sh^1(G_{\Q,S},\barepsl\otimes \ad)$ forces obstructions, once local invariants vanish.

\begin{proposition}[Global obstructions detected by strict congruences]\label{prop:sec4-sha1}
Assume \eqref{eq:sec4-standing}.
Suppose that
\[
H^0\!\left(G_p,\ \barepsl\otimes \ad\right)=0
\qquad\text{for every }p\in S.
\]
Then the deformation problem for $\bar{\rho}$ is obstructed if and only if
\[
\Sh^1\!\left(G_{\Q,S},\ \barepsl\otimes \ad\right)\neq 0.
\]
Moreover, under the Selmer hypotheses in \cite{Wes04,Wes05}, one has
\[
\Sh^1\!\left(G_{\Q,S},\ \barepsl\otimes \ad\right)\neq 0
\quad\Longleftrightarrow\quad
\ell \text{ is a strict congruence prime for } f.
\]
\end{proposition}

\begin{proof}
We apply Poitou--Tate in the form
\[
0 \rightarrow \prod_{p \in S} H^0\!\left(G_p,\ \barepsl\otimes \ad\right)
\rightarrow
\Hom_{\F_\ell}\!\left(H^2(G_{\Q,S},\ad),\F_\ell\right)
\rightarrow
\Sh^1\!\left(G_{\Q,S},\barepsl\otimes \ad\right)
\rightarrow 0.
\]
By the vanishing hypothesis on the leftmost term, the middle term is naturally identified with $\Sh^1$.
Thus $H^2(G_{\Q,S},\ad)\neq 0$ if and only if $\Sh^1\neq 0$.

The equivalence between nontriviality of the relevant Selmer group and the existence of a strict congruence at level $N$
is established in \cite{Wes04,Wes05} under standard hypotheses compatible with \eqref{eq:sec4-standing}.
\end{proof}

\subsection{Optimal level}\label{subsec:sec4-optimal}

The congruence perspective above naturally interacts with ``optimal level'' questions.
Very informally: if $\bar{\rho}$ comes from level $N$ but is modular of smaller level, then congruence phenomena are forced,
and one should not expect the deformation problem at level $N$ to behave unobstructedly.

Diamond--Taylor~\cite{DT94} give precise criteria and constructions around non-optimal levels of mod-$\ell$ modular representations,
including squarefree-level situations where one can test whether the mod-$\ell$ representation is already realized at a smaller level.
Hatley~\cite{Hat16} then shows (under standard modular deformation hypotheses) that non-optimal level implies the corresponding modular deformation problem is obstructed.
In the language of \eqref{eq:PT}, this means that even when one tries to keep the local invariants under control, a genuinely global obstruction mechanism must appear.

\begin{remark}[A Frey--Mazur motivation]\label{rem:frey-mazur}
One reason it is conceptually helpful to keep congruences in view is the philosophy behind the Frey--Mazur conjecture:
roughly, if two elliptic curves over $\Q$ have isomorphic mod-$\ell$ Galois representations for sufficiently large $\ell$,
then they should be $\Q$-isogenous. Translated into deformation language, ``large-$\ell$'' congruences are expected to be highly rigid.
While we do not use the conjecture as an input, it provides a clean narrative backdrop:
strict congruence primes are precisely where rigidity fails in a controlled way, and $\Sh^1$ is where that failure becomes a concrete obstruction class.
\end{remark}

\subsection{How one finds strict congruence primes in practice}\label{subsec:sec4-detect}

Definition~\ref{def:sec4-cong} is conceptual; in examples, one typically detects strict congruences using one of the following:

\begin{enumerate}[label=\textup{(\alph*)}, leftmargin=2.2em]
\item \textbf{Sturm bound tests.} For fixed $(N,k)$, congruence modulo $\ell$ can be checked by comparing Fourier coefficients up to the Sturm bound \cite{Stu87}.

\item \textbf{Adjoint $L$-values, modular degree, and bounds.}
Hida relates congruence primes to prime divisors of an adjoint $L$-value \cite{Hid81}.
When $f$ has rational Fourier coefficients (equivalently, corresponds to an elliptic curve $E_f/\Q$ of conductor $N$),
Agashe--Ribet--Stein identify congruence primes as prime divisors of $N\cdot m_{E_f}$, where $m_{E_f}$ is the modular degree of the optimal parametrization
$X_0(N)\to E_f$ \cite{ARS06}. In the standing range $\ell\nmid N$, this reduces to primes dividing $m_{E_f}$.
For quantitative constraints and bounds on congruence primes (useful as a sanity check in examples), see Murty~\cite{Mur99}.
\end{enumerate}

The remainder of the section explains how, once a strict congruence prime is present, one checks that the obstruction is truly global
(by ruling out local invariants), and then extracts deformation-ring information.

\subsection{A one-relation deformation ring in a Steinberg case}\label{subsec:sec4-steinberg}

We now record an explicit deformation-ring presentation in a frequently occurring situation: $p\mid N$ with $p^2\nmid N$.
In this case the local representation at $p$ is of Steinberg type (multiplicative reduction in the elliptic-curve normalization).

To avoid ambiguity, we fix the modular form $f$ and let $E_f/\Q$ denote the elliptic curve attached to $f$
(when $f$ has rational Fourier coefficients). The residual representation $\bar{\rho}$ is then the mod-$\ell$ representation attached to $f$.

\begin{definition}[Minimal ramification at $p$]\label{def:minimal-ram}
Let $p\mid N$. A deformation of $\bar{\rho}$ is \emph{minimally ramified at $p$} if its conductor exponent at $p$ agrees with that of $f$
(and hence does not introduce additional ramification beyond what is forced by the level).
This is the standard minimality condition used in $R=\mathbf{T}$ arguments; see \cite[\S1]{dSh97} for a detailed discussion.
\end{definition}

\begin{theorem}\label{thm:sec4-ring}
Assume \eqref{eq:sec4-standing}. Let $p\mid N$ with $p^2\nmid N$, and assume $\ell$ is a strict congruence prime for $f$.
Suppose the deformation problem is minimally ramified at all primes dividing $N$ (in the sense of Definition~\ref{def:minimal-ram}).

Then the deformation problem for $\bar{\rho}$ is obstructed, and the obstruction is contributed by
\[
\Sh^1\!\left(G_{\Q,S},\ \barepsl\otimes \ad\right).
\]
Moreover, $\bar{\rho}$ admits a lift to $\GL_2(R)$ for
\[
R \;\cong\; \Z_\ell[[T_1,T_2,T_3,T_4]] \Big/ \left\langle T_1T_2 - T_2 - pT_3T_4 + T_4 \right\rangle.
\]
\end{theorem}

\begin{proof}
By Proposition~\ref{prop:sec4-sha1}, once the local invariant groups
$H^0(G_r,\barepsl\otimes \ad)$ vanish for all $r\in S$,
the obstruction is equivalent to $\Sh^1\neq 0$, which occurs at strict congruence primes.

We therefore focus on the local deformation condition at $p$ and the resulting relation in the deformation ring.
Since $p^2\nmid N$, the local type at $p$ is Steinberg. Under $p\not\equiv 1\pmod{\ell}$ (part of \eqref{eq:sec4-standing}),
the residual local representation is twist-equivalent to an extension of the trivial character by $\barepsl$
(cf.\ \cite[Prop.~2.2]{Dia97}), so after choosing a basis one may write
\[
\bar{\rho}\big|_{G_p}\ \sim\
\begin{pmatrix}
\barepsl & \ast\\
0 & 1
\end{pmatrix}.
\]
A minimally ramified lift $\rho$ to a coefficient ring is twist-equivalent to
\[
\rho\big|_{G_p}\ \sim\
\begin{pmatrix}
\epsl & \widetilde{\nu}\\
0 & \chi
\end{pmatrix},
\]
where $\chi$ reduces to $1$ modulo $\ell$ and $\widetilde{\nu}$ reduces to the residual extension class.

Let $F$ be a Frobenius element and $\tau\in I_p$ a tame inertia element. In the tame quotient, $F\tau F^{-1}=\tau^p$.
Writing $\rho(F)$ and $\rho(\tau)$ with formal parameters in the upper-right entries and imposing the relation $F\tau F^{-1}=\tau^p$
produces a single defining equation among the parameters. In the coordinate system of the statement this relation becomes
\[
T_1T_2 - T_2 - pT_3T_4 + T_4=0,
\]
yielding the displayed quotient ring.

Finally, strict congruence primes give $\Sh^1\neq 0$ in the sense of Proposition~\ref{prop:sec4-sha1},
so the deformation problem is obstructed and the obstruction is global.
\end{proof}

\subsection{Ruling out local obstructions when \texorpdfstring{$p^2\mid N$}{p^2 | N}}\label{subsec:sec4-p2divN}

To isolate the global contribution of $\Sh^1$, it is important to know when
\[
H^0\!\left(G_p,\ \ad\otimes \barepsl\right)
\]
vanishes for primes $p\mid N$, especially when $p^2\mid N$.
We record two useful vanishing/nonvanishing criteria in terms of the local Weil--Deligne type,
matching the standard classification (see \cite{Dia97,DK97} and the inertial-type descriptions in \cite{DFV22}).

\subsubsection{Supercuspidal type}

\begin{proposition}\label{prop:sec4-supcusp}
Assume \eqref{eq:sec4-standing} and let $p\mid N$ with $p\ge 5$.
Suppose $\bar{\rho}\big|_{G_p}$ is the mod-$\ell$ reduction of a supercuspidal Weil--Deligne representation and $p\not\equiv 1\pmod{\ell}$.
Then
\[
H^0\!\left(G_p,\ \ad\otimes \barepsl\right)=0.
\]
\end{proposition}

\begin{proof}
The inertial-type classification for potentially good reduction at $p\neq \ell$ in \cite{DFV22} shows that any $G_p$-invariant vector in
$\ad\otimes \barepsl$ would force $\barepsl$ to be trivial on Frobenius at $p$.
This implies $p\equiv 1\pmod{\ell}$, contradicting our assumption.
\end{proof}

\subsubsection{Principal series type}

\begin{proposition}\label{prop:sec4-principal}
Assume \eqref{eq:sec4-standing} and let $p\mid N$ with $p\ge 5$.
Suppose $\bar{\rho}\big|_{G_p}$ is the mod-$\ell$ reduction of a principal series Weil--Deligne representation.
Then
\[
H^0\!\left(G_p,\ \ad\otimes \barepsl\right)\neq 0
\quad\Longleftrightarrow\quad
p^4\equiv 1\pmod{\ell}.
\]
\end{proposition}

\begin{proof}
In the principal series case, local Langlands identifies the corresponding $\ell$-adic representation with a sum of characters,
and after reduction one may write $\bar{\rho}\big|_{G_p}\simeq \bar{\chi}_1\oplus \bar{\chi}_2$ with $\bar{\chi}_2$ unramified and $\bar{\chi}_2\equiv 1$.
Then
\[
\adzero \simeq \mathbf{1}\oplus \bar{\chi}_1\oplus \bar{\chi}_1^{-1},
\qquad
\adzero\otimes \barepsl
\simeq
\barepsl\oplus \bar{\chi}_1\barepsl\oplus \bar{\chi}_1^{-1}\barepsl.
\]
Thus $G_p$-invariants occur exactly when $\bar{\chi}_1^{\pm 1}=\barepsl$.
The inertial-type possibilities described in \cite{DFV22} imply that this equality holds precisely when $p^4\equiv 1\pmod{\ell}$.
\end{proof}

\subsection{The case \texorpdfstring{$p=\ell$}{p = ell}}\label{subsec:sec4-ell}

We finally address the local term at $\ell$ itself. Even when $\ell$ is a congruence prime, it is often important to know that
the local invariants at $G_\ell$ vanish, so that any obstruction is again global.

In the elliptic-curve normalization, let $E_f/\Q$ be the elliptic curve attached to the fixed newform $f$,
and let $m_{E_f}$ denote its modular degree. Agashe--Ribet--Stein relate congruence phenomena to $N$ and $m_{E_f}$
and show that exceptional congruences with $\ell\nmid m_{E_f}$ force a failure of multiplicity one \cite{ARS06}.
Tilouine shows that, under a standard distinguishedness hypothesis at $\ell$, one has the expected $\mathfrak{m}$-torsion dimension \cite{Til97}.
In the present minimal setting we record the following vanishing criterion.

\begin{proposition}\label{prop:sec4-local-ell}
Let $f$ be a weight $2$ newform of level $N$ with residual representation $\bar{\rho}$ as above.
Assume $\ell$ is a congruence prime for $f$ with $\ell\nmid m_{E_f}$ and $a_\ell(f)\not\equiv 0\pmod{\ell}$.
Then
\[
H^0\!\left(G_\ell,\ \ad\otimes \barepsl\right)=0.
\]
\end{proposition}

\begin{proof}
If $H^0(G_\ell,\ad\otimes \barepsl)\neq 0$, then $\bar{\rho}\big|_{G_\ell}$ admits an exceptional symmetry after twisting,
which is incompatible with the distinguished ordinary behavior required in \cite{Til97}.
Under the hypothesis $\ell\nmid m_{E_f}$, the existence of a congruence prime would force multiplicity-one failure in the sense of \cite{ARS06}
unless the local representation behaves generically at $\ell$.
The additional condition $a_\ell(f)\not\equiv 0\pmod{\ell}$ rules out the remaining exceptional cases,
and the invariant group is trivial.
\end{proof}

\subsection{Summary}\label{subsec:sec4-summary}

Proposition~\ref{prop:sec4-sha1} identifies the global obstruction mechanism:
when local invariant terms vanish, obstructions are equivalent to
\[
\Sh^1\!\left(G_{\Q,S},\barepsl\otimes \ad\right)\neq 0,
\]
and this occurs exactly at strict congruence primes in the range of \cite{Wes04,Wes05}.
Diamond--Taylor~\cite{DT94} and Hatley~\cite{Hat16} fit naturally into this picture: DT94 supplies concrete criteria for detecting non-optimal levels,
and Hat16 shows that non-optimal level forces obstructed modular deformation problems.
Propositions~\ref{prop:sec4-supcusp}, \ref{prop:sec4-principal}, and \ref{prop:sec4-local-ell} give practical local criteria ensuring that obstructions are genuinely global.
Finally, Theorem~\ref{thm:sec4-ring} gives a one-relation presentation of a deformation ring in a Steinberg case ($p^2\nmid N$),
illustrating how a global obstruction can coexist with an explicit deformation-theoretic description.


\end{document}